\documentclass[11pt]{amsart}
\usepackage{amsmath, amssymb, latexsym, amsthm, mathrsfs}

\newcommand{\ed}{

\end{document}
}

\newcommand{\Arh}{Arhangel'ski\u{\i}}

\newcommand{\Setting}[7]{\xymatrix@R=4pt@C=7pt{#1\ar@{-}[r]&#2\ar@{-}[r]&#3\\&#4\ar@{-}[u]\\
#5\ar@{-}[uu]\ar@{-}[r] & #6\ar@{-}[u]\ar@{-}[r] & #7\ar@{-}[uu]}}

\newcommand{\EF}{\mathbb{EF}}

\newcommand{\compactN}{\cl{\mathbb{N}}}

\newcommand{\bq}{\begin{quote}}
\newcommand{\eq}{\end{quote}}
\newcommand{\cl}[1]{\overline{#1}}

\newcommand{\inv}{^{-1}}

\newcommand{\N}{\mathbb{N}}
\newcommand{\NN}{{\N^{\N}}}

\newcommand{\NcompactN}{{\compactN^\N}}

\newcommand{\sseq}[1]{\{#1 : n\in\N\}}

\newcommand{\cI}{\mathcal{I}}

\newcommand{\scrA}{\mathscr{A}}
\newcommand{\scrB}{\mathscr{B}}

\newcommand{\rmB}{\mathrm{B}}
\newcommand{\rmC}{\mathrm{C}}
\newcommand{\rmF}{\mathrm{F}}

\newcommand{\BG}{\rmB_\Gamma}
\newcommand{\FG}{\rmF_\Gamma}

\newcommand{\CG}{\rmC_\Gamma}

\newcommand{\cF}{\mathcal{F}}

\newcommand{\rmO}{\mathrm{O}}

\newcommand{\R}{\mathbb{R}}
\newcommand{\cU}{\mathcal{U}}
\newcommand{\cC}{\mathcal{C}}
\newcommand{\Union}{\bigcup}
\newcommand{\cV}{\mathcal{V}}

\newcommand{\Impl}{\Rightarrow}
\long\def\forget#1\forgotten{}

\newcommand{\oo}{\infty}

\newcommand{\x}{\times}
\newcommand{\Iff}{\Leftrightarrow}
\newcommand\comp{^{\text{\tt c}}}
\newcommand{\nin}{\notin}

\newcommand{\sbst}{\subseteq}
\newcommand{\spst}{\supseteq}
\newcommand{\sm}{\setminus}
\newcommand{\as}{\subseteq^*}

\newtheorem{thm}{Theorem}
\newcommand{\bthm}{\begin{thm}} \newcommand{\ethm}{\end{thm}}
\newtheorem{prop}[thm]{Proposition}
\newcommand{\bprp}{\begin{prop}} \newcommand{\eprp}{\end{prop}}
\newtheorem{fact}[thm]{Fact}
\newcommand{\bfct}{\begin{fact}} \newcommand{\efct}{\end{fact}}
\newtheorem{prob}[thm]{Problem}
\newcommand{\bprb}{\begin{prob}} \newcommand{\eprb}{\end{prob}}
\newtheorem{lem}[thm]{Lemma}
\newcommand{\blem}{\begin{lem}} \newcommand{\elem}{\end{lem}}
\newtheorem{cor}[thm]{Corollary}
\newcommand{\bcor}{\begin{cor}} \newcommand{\ecor}{\end{cor}}
\newtheorem{conj}[thm]{Conjecture}
\newcommand{\bcnj}{\begin{conj}} \newcommand{\ecnj}{\end{conj}}
\theoremstyle{definition}
\newtheorem{defn}[thm]{Definition}
\newcommand{\bdfn}{\begin{defn}} \newcommand{\edfn}{\end{defn}}
\theoremstyle{remark}
\newtheorem{rem}[thm]{Remark}
\newcommand{\brem}{\begin{rem}} \newcommand{\erem}{\end{rem}}
\newtheorem{exam}[thm]{Example}
\newcommand{\bexs}{\begin{exam}} \newcommand{\eexs}{\end{exam}}
\newcommand{\bpf}{\begin{proof}} \newcommand{\epf}{\end{proof}}

\newcommand{\be}{\begin{enumerate}}
\newcommand{\ee}{\end{enumerate}}
\newcommand{\bi}{\begin{itemize}}
\newcommand{\itm}{\item}
\newcommand{\ei}{\end{itemize}}

\newcommand{\sone}{\mathsf{S}_1}

\newcommand{\ufin}{\mathsf{U}_\mathrm{fin}}

\title[Hurewicz and sheaf amalgamations]{Hereditarily Hurewicz spaces and Arhangel'ski\u{\i} sheaf amalgamations}

\author{Boaz Tsaban}
\address[Tsaban]{Department of Mathematics,
Bar-Ilan University, Ramat Gan 52900, Israel;
and Department of Mathematics,
Weizmann Institute of Science, Rehovot 76100, Israel.}
\email{tsaban@math.biu.ac.il}
\thanks{Partially supported by the Koshland Center for Basic Research.}

\author{Lyubomyr Zdomskyy}
\address[Zdomskyy]{Department of Mathematics,
Weizmann Institute of Science, Rehovot 76100, Israel.}
\curraddr{Universit\"at Wien,
Kurt G\"odel Research Center for Mathematical Logic, W\"ahringer Str.\ 25, A-1090 Wien,
Austria.}
\email{lzdomsky@gmail.com}



\begin{document}

\begin{abstract}
A classical theorem of Hurewicz characterizes spaces with the
Hurewicz covering property as those having bounded continuous
images in the Baire space. We give a similar characterization for
spaces $X$ which have the Hurewicz property hereditarily.

We proceed to consider the class of \Arh{} $\alpha_1$ spaces, for
which every sheaf at a point can be amalgamated in a natural way.
Let $C_p(X)$ denote the space of continuous real-valued functions
on $X$ with the topology of pointwise convergence. Our main result
is that $C_p(X)$ is an $\alpha_1$ space if, and only if, each
\emph{Borel} image of $X$ in the Baire space is bounded. Using
this characterization, we solve a variety of problems posed in the
literature concerning spaces of continuous functions.
\end{abstract}

\maketitle

\tableofcontents

\section{Introduction}

We are mainly concerned with spaces $X$ which are (homeomorphic
to) sets of irrational numbers, and we recommend adopting this
restriction for clarity. Our results (and proofs) apply to all
topological spaces $X$ in which each open set is a union of
countably many clopen sets, and the spaces considered are assumed
to have this property.\footnote{Every perfectly
normal space (open sets are $F_\sigma$) with upper inductive
dimension $0$ (disjoint closed sets can be separated by a clopen
set) has the required property. Thus, the spaces considered in the
references also have the required property.}

Fix a topological space $X$. Let $\scrA,\scrB$ be families of covers of $X$.
The space $X$ may or may not have the following property \cite{coc1}.
\begin{description}
\itm[$\ufin(\scrA,\scrB)$]
Whenever $\cU_1,\cU_2,\dots\in\scrA$ and none contains a
finite subcover, there exist finite sets $\cF_n\sbst\cU_n$, $n\in\N$,
such that $\sseq{\Union\cF_n}\in\scrB$.
\end{description}
Let $\rmO$ denote the collection of all countable open covers of $X$.\footnote{If
$X$ is Lindel\"of, we can consider arbitrary open covers of $X$.}
A cover $\cU$ of $X$ is \emph{point-cofinite} if $\cU$ is infinite and each $x\in X$,
is a member of all but finitely many members of $\cU$.\footnote{Traditionally, point-cofinite
covers were called \emph{$\gamma$-covers} \cite{GN}.}
Let $\Gamma$ denote the collection of all open point-cofinite covers of $X$.
Motivated by studies of Menger \cite{Menger24}, Hurewicz \cite{Hure27}
introduced the \emph{Hurewicz property} $\ufin(\rmO,\Gamma)$.

Hurewicz \cite{Hure27} essentially obtained the following combinatorial
characterization of $\ufin(\rmO,\Gamma)$ (see Rec\l{}aw \cite{Rec94}).
For $f,g\in\NN$, $f\le^* g$ means $f(n)\le g(n)$ for all but
finitely many $n$. A subset $Y$ of $\NN$ is \emph{bounded} if
there is $g\in\NN$ such that $f\le^* g$ for all $f\in Y$.

\bthm[Hurewicz]\label{hure}
$X$ satisfies $\ufin(\rmO,\Gamma)$ if, and only if, every
continuous image of $X$ in $\NN$ is bounded.
\ethm
This characterization has found numerous applications---see
\cite{LecceSurvey, KocSurv, ict} and references therein.
We give a similar characterization for \emph{hereditarily Hurewicz} spaces,
that is, spaces $X$ such that each subspace of $X$ satisfies
$\ufin(\rmO,\Gamma)$.

The property of being hereditarily Hurewicz was studied in, e.g.,
\cite{FM88, Nowik99, Miller07}.
Rubin introduced a property of subsets of $\R$
such that the existence of a set with this property
is equivalent to the possibility of a certain construction of boolean algebras \cite{Rubin83}.
Miller \cite{Miller07} proved that the \emph{Rubin spaces}
are exactly the hereditarily Hurewicz spaces.

The property of being hereditarily Hurewicz also manifests itself
as follows:
A set $X\sbst\R$ is a \emph{$\sigma'$ space}\label{sig'def}
\cite{Sakai98} if
for each $F_\sigma$ set $E$, there is an $F_\sigma$ set $F$
such that $E\cap F=\emptyset$ and $X\sbst E\cup F$.
This property was effectively used in studies of generalized metric spaces
\cite{Chen03}.
Recently, Sakai proved that $X$ is a $\sigma'$ space if, and only if,
$X$ is hereditarily Hurewicz (Theorem \ref{sig'} below).

There exist additional classes of hereditarily Hurewicz spaces in
the literature. We describe some of them.

A topological space is \emph{Fr\'echet} if each point in the closure of a subset
of the space is a limit of a convergent sequence of points from that subset.
The following concepts, due to \Arh{} \cite{Arh72, Arh79}, are important in determining
when a product of Fr\'echet spaces is Fr\'echet.
Let $Y$ be a general topological space (not necessarily Lindel\"of or zero-dimensional).
A \emph{sheaf} at a point $y\in Y$ is a family of sequences, each converging to $y$.
To avoid trivialities, we consider only sequences of distinct elements.
We say that a countable set $A$ converges to $y$ if some (equivalently, each) bijective
enumeration of $A$ converges to $y$.
The space $Y$ is an \emph{$\alpha_1$ space} if for each $y\in Y$, each countable sheaf $\sseq{A_n}$ at $y$
can be amalgamated as follows: There are cofinite subsets $B_n\sbst A_n$, $n\in\N$, such that the set
$B=\Union_nB_n$ converges to $y$.
The references dealing with $\alpha_1$ spaces are too numerous to be listed here;
see \cite{Shakhmatov02} and the references therein for a partial list.

Fix a space $X$. The space $C_p(X)$ is the family of all continuous real-valued functions on $X$,
viewed as a subspace of the Tychonoff product $\R^X$.
A sequence of results by Bukovsk\'y--Rec\l{}aw--Repick\'y \cite{BRR91}, Rec\l{}aw \cite{Rec97},
Sakai \cite{Sakai07}, and Bukovsk\'y--Hale\v{s} \cite{BH07}, culminated in the result
that if $C_p(X)$ is an $\alpha_1$ space, then $X$ is hereditarily Hurewicz.
Our main result is that if $C_p(X)$ is an $\alpha_1$ space, then each Borel image
of $X$ in $\NN$ is bounded. It is easy to see that the converse implication also holds,
and we obtain a powerful characterization of spaces $X$ such that $C_p(X)$ is an $\alpha_1$ space.

Historically, the realization that if $C_p(X)$ is an $\alpha_1$ space then $X$ is hereditarily Hurewicz
goes through QN spaces \cite{BRR91}:\label{QNdef}
Let $Y$ be a metric space.
A function $f:X\to Y$ is a \emph{quasi-normal} limit of functions $f_n:X\to Y$
if there are positive reals $\epsilon_n$, $n\in\N$, converging to $0$
such that for each $x\in X$, $d(f_n(x),f(x))<\epsilon_n$ for all but finitely
many $n$.
A topological space $X$ is a \emph{QN space} if whenever $0$ is a pointwise limit
of a sequence of continuous real-valued functions on $X$,
we have that $0$ is a quasi-normal limit of the same sequence.
QN spaces are studied in, e.g., \cite{BRR91, Rec97, alpha_i, Nowik99, BRR01, Sakai07, BH07}.
In \cite{Sakai07, BH07} it was shown that $X$ is a QN space if, and only if, $C_p(X)$ is an $\alpha_1$ space.
Thus, QN spaces are also characterized by having bounded Borel images in $\NN$.

We use our main theorem to show that quite a few additional properties studied in the literature
are equivalent to having bounded Borel images in $\NN$, and consequently
solve a variety of problems posed in the literature.
To make the paper self-contained and accessible to a wide audience,
we supply proofs for all needed results. Often, our proofs of known results are slightly simpler
than those available in the literature.

\section{A characterization of hereditarily Hurewicz spaces}

Let $\compactN=\N\cup\{\oo\}$ be the one-point compactification
of $\N$, and endow $\NcompactN$ with the Tychonoff product topology.
An element $f\in\NcompactN$ is \emph{eventually finite} if there is $m$
such that $f(n)<\oo$ for all $n\ge m$.
Let $\EF$ be the subspace of $\NcompactN$ consisting of all eventually finite
elements of $\NcompactN$. The partial order $\le^*$ extends to $\EF$ in the natural way.

\bthm\label{main}
$X$ is hereditarily $\ufin(\rmO,\Gamma)$ if, and only if, every
continuous image of $X$ in $\EF$ is bounded.
\ethm
\bpf
$(\Impl)$ Assume that $\Psi:X\to\EF$ is continuous.
For each $n$, define the following ($G_\delta$) subset of $\Psi[X]$:
$$G_n = \{f\in\Psi[X] : (\forall m\ge n)\ f(m)<\oo\}.$$
Let $T^n:G_n\to\NN$ be the shift transformation
defined by $T^n(f)(m)=f(m+n)$ for all $m$.

For each $n$, $X_n=\Psi\inv[G_n]\sbst X$,
and therefore satisfies
$\ufin(\rmO,\Gamma)$. By Theorem \ref{hure},
$T^n[\Psi[X_n]]=T^n[G_n]$ is a bounded subset of $\NN$.
Thus, $G_n$ is a bounded subset of $\EF$, and therefore so is
$\Psi[X]=\Union_n G_n$.

$(\Leftarrow)$ First, note that Hurewicz's Theorem \ref{hure} and our assumption on $X$ imply that
$X$ satisfies $\ufin(\rmO,\Gamma)$.

\blem\label{Gdelta}
If each $G_\delta$ subset of $X$ satisfies $\ufin(\rmO,\Gamma)$, then $X$ is hereditarily $\ufin(\rmO,\Gamma)$.
\elem
\bpf
Let $Y\sbst X$.
Assume that $\cU_n$, $n\in\N$, are covers of $Y$ by open subsets of $X$,
which do not contain finite subcovers.
Each $\cU_n$ is an open cover of $G=\bigcap_n\Union\cU_n\spst Y$, and has no
finite subcover of $G$.
As $G$ is a $G_\delta$ subset of $X$, it satisfies $\ufin(\rmO,\Gamma)$.
Thus, there are finite $\cF_n\sbst\cU_n$, $n\in\N$,
such that $\sseq{\Union\cF_n}$ is a point-cofinite cover of $G$, and therefore of $Y$.
\epf

Assume that $G$ is a $G_\delta$ subset of $X$.

\blem[Sakai \cite{Sakai07}]\label{SakaiGdelta}
For each $G_\delta$ subset $G$ of $X$, there is an open point-cofinite cover $\sseq{U_n}$ of $X$ such that $G=\bigcap_nU_n$.
\elem
\bpf
$G^c=\Union_n C_n$ with each $C_n$ closed.
If $A$ is closed and $B$ is open, then $B$ is a union of countably many disjoint clopen sets,
and therefore $A\cap B$ is a union of countably many disjoint closed sets.
Thus, each of the disjoint sets $C_n\sm(C_1\cup\dots,C_{n-1})$, $n\in\N$,
is a union of countably many disjoint closed sets. Hence, $G^c=\Union_n\tilde C_n$ where the sets
$\tilde C_n$ are closed and disjoint, and therefore
$G=\bigcap_n\tilde C_n^c$, where $\sseq{\tilde C_n^c}$ is
an open point-cofinite cover of $X$.
\epf

So, let $\sseq{U_n}$ be an open point-cofinite cover of $X$ such that $G=\bigcap_nU_n$.
For each $n$, let $U_n=\Union_m C^n_m$, a union of disjoint clopen sets.
Define $\Psi:X\to\EF$ by
$$\Psi(x)(n) = \begin{cases}
m & m\in\N, x\in C^n_m\\
\oo & x\nin U_n
\end{cases}
$$
As $\sseq{U_n}$ is a point-cofinite cover of $X$, $\Psi(x)$ is eventually
finite for each $x\in X$.

The function $\Psi$ is continuous:
A basic open set in $\EF$ has the form $\prod_n V_n$ such
that there are finite $I_0,I_1\sbst\N$ and elements $m_n$, $n\in I_0\cup I_1$,
for which: For each $n\in I_0$, $V_n=\{m_n\}$, for each $n\in I_1$, $V_n=\{m_n, m_n+1,\dots\}\cup\{\oo\}$,
and for each $n\nin I_0\cup I_1$, $V_n=\N\cup\{\oo\}$.
Now,
$$\Psi\inv\left [\prod_{n\in\N} V_n\right ] =
\bigcap_{n\in I_0} C^n_{m_n}\cap\bigcap_{n\in I_1}\left(X\sm\left(\Union_{k<m_n}C^n_k\right)\right)$$
is open.

Thus, $\Psi[X]$ is bounded by some $g\in\NN$.
Now,
$$G = \{x\in X : (\forall n)\ \Psi(x)(n)<\oo\} = \Psi\inv[\{f\in\NN : f\le^* g\}].$$
The set $\{f\in\NN : f\le^* g\}$ is an $F_\sigma$ subset of $\EF$. Indeed, let
$\sseq{g_n}$ enumerate all elements of $\NN$ which are eventually equal to $g$.
Then $\{f\in\NN : f\le^* g\} = \Union_n \{f\in\EF : f\le g_n\}$.
Thus, $G$ is an $F_\sigma$ subset of $X$.
As $\ufin(\rmO,\Gamma)$ is hereditary for closed subsets and preserved by countable unions,
$G$ satisfies $\ufin(\rmO,\Gamma)$.
\epf

Recall that a topological space $X$ is a \emph{$\sigma$ space} if each
$G_\delta$ subset of $X$ is an $F_\sigma$ subset of $X$.
The proof of Theorem \ref{main} actually shows that
$(2\Impl 3)$, $(2\Impl 1)$, $(3\Impl 4)$, and $(4\Impl 2)$
in the following theorem (and therefore establishes it).

\bthm\label{manyequivs}
The following are equivalent:
\be
\itm $X$ is hereditarily $\ufin(\rmO,\Gamma)$.
\itm Each $G_\delta$ subset of $X$ satisfies $\ufin(\rmO,\Gamma)$.
\itm Every continuous image of $X$ in $\EF$ is bounded.
\itm $X$ satisfies $\ufin(\rmO,\Gamma)$ and is a $\sigma$ space.
\hfill\qed
\ee
\ethm
The implication $(1\Impl 4)$ in Theorem \ref{manyequivs} was first proved
by Fremlin and Miller \cite{FM88}. The implication $(4\Impl 1)$ can be
alternatively deduced from Theorem 3.12 of \cite{BRR01} and Corollary 10 of \cite{BH03}.
An additional equivalent formulation was discovered by Sakai.
Recall the definition of $\sigma'$ space from the introduction (page \pageref{sig'def}).

\bthm[Sakai]\label{sig'}
Let $X\sbst\R$.
$X$ is a $\sigma'$ space if, and only if, $X$ is hereditarily $\ufin(\rmO,\Gamma)$.
\ethm
\bpf
This follows from Theorem 5.7 of \cite{coc2}: $X$ satisfies $\ufin(\rmO,\Gamma)$
if, and only if, for each $G_\delta$ set $G\sbst\R$ containing $X$, there is an $F_\sigma$ set
$F\sbst\R$ such that $X\sbst F\sbst G$.

$(\Impl)$ As being a $\sigma'$ space is hereditary, it suffices to show that
$X$ satisfies $\ufin(\rmO,\Gamma)$. Indeed, for each $G_\delta$ set $G\sbst\R$ containing $X$,
let $E=\R\sm G$, and take an $F_\sigma$ set $F\sbst\R$ disjoint from $E$ such that
$X\sbst E\cup F$. Then $X\sbst F\sbst G$.

$(\Leftarrow)$ Let $E\sbst\R$ be $F_\sigma$.
As $X\sm E$ satisfies $\ufin(\rmO,\Gamma)$ and is a subset of the $G_\delta$ set $\R\sm E$,
there is an $F_\sigma$ set $F\sbst\R$ such that $X\sm E\sbst F\sbst\R\sm E$.
Then $E\cap F=\emptyset$ and $X\sbst E\cup F$.
\epf

To indicate the potential usefulness of Theorem \ref{main},
we use it to give slightly more direct proofs of two known theorems.
Recall the definition of QN spaces from the introduction (page \pageref{QNdef}).

\bthm[Rec\l{}aw \cite{Rec97}]\label{RecHered}
If $X$ is a QN space, then $X$ is hereditarily $\ufin(\rmO,\Gamma)$.
\ethm
\bpf
Let $Y\sbst\EF$ be a continuous image of $X$.
Then $Y$ is a QN space. By Theorem \ref{main}, it suffices to show that $Y$ is bounded.
For each $n$, and each $y\in Y$, define
$$f_n(y) = \frac{1}{\min y\inv(n)},$$
using the natural conventions that $\min\emptyset=\oo$ and $1/\oo=0$.
$\lim_n f_n(y)\allowbreak=0$ for all $y\in Y$.
As $Y$ is a QN space,
there are positive $\epsilon_n$, $n\in\N$, dominating this convergence.
For each $k$, let
$$Y_k=\{y\in Y : (\forall n\ge k)\ f_n(y)<\epsilon_n\}.$$
$Y=\Union_k Y_k$. We will show that each $Y_k$ is bounded.

Fix $k$. Take an increasing $g\in\NN$ such that $g(1)=k$ and for each $n$,
$\epsilon_m<1/n$ for all $m\ge g(n)$. Then $Y_k$ is bounded by $g$:
Let $y\in Y_k$. Fix $n$ such that $y(n)<\oo$. If $y(n)\le k$, then
$y(n)\le g(1)\le g(n)$. Otherwise, $y(n)>k$, and since $y\in Y_k$,
$f_{y(n)}(y)<\epsilon_{y(n)}$. Thus,
$$\frac{1}{n}\le \frac{1}{\min y\inv (y(n))}=f_{y(n)}(y)<\epsilon_{y(n)},$$
and therefore $y(n)$ cannot be greater than $g(n)$.
\epf

\bthm[Rec\l{}aw \cite{Rec97}]\label{reclaw}
If $X$ is QN space, then $X$ is a $\sigma$ space.
\ethm
\bpf
Theorems \ref{manyequivs} and \ref{RecHered}.
\epf

The following sections give a deeper reason for the last two theorems.

\section{Bounded Borel images}

Our main goal in this section is to establish the equivalence in the following Theorem \ref{Borel}.
The implication $(2\Impl 1)$ in this theorem is Proposition 9 of Scheepers \cite{alpha_i}.
The implication $(1\Impl 2)$ is the more difficult one, and will
be proved in a sequence of related results.

\bthm\label{Borel}
The following are equivalent:
\be
\itm $C_p(X)$ is an $\alpha_1$ space.
\itm Each Borel image of $X$ in $\NN$ is bounded.
\ee
\ethm
\bpf
$(2\Impl 1)$ Consider a sheaf $\sseq{A_n}$ at $f\in C_p(X)$.
For each $n$, enumerate $A_n = \{f^n_m : m\in\N\}$ bijectively.
Define a Borel function $\Psi:X\to\NN$ by
$$\Psi(x)(n) = \min\{k : (\forall m\ge k)\ |f^n_m(x)-f(x)|\le 1/n\}.$$
Let $g\in\NN$ bound $\Psi[X]$, and take the amalgamation
$B=\Union_n\{f^n_m : m\ge g(n)\}$. Then $B$ converges to $f$.

$(1\Impl 2)$ Assume that $C_p(X)$ is an $\alpha_1$ space.
Then the subspace $C_p(X,\{0,1\})$ of $C_p(X)$,
consisting of all continuous functions $f:X\to\{0,1\}$, is an $\alpha_1$ space.\footnote{In fact,
by the methods of Gerlits--Nagy \cite{GNFr}, the converse implication also holds. This fact will not
be used in our proof.}
Each element of $C_p(X,\{0,1\})$ has the form $\chi_U$, the characteristic function of
a clopen set $U\sbst X$.
Immediately from the definition, a sequence $\chi_{U_n}$ of elements of
$C_p(X,\{0,1\})$ converges pointwise to the constant function
$1$ if, and only if, $\sseq{U_n}$ is a clopen point-cofinite cover of $X$.
This gives the following, which is due to Bukovsk\'y--Hale\v{s} (cf. \cite[Theorem 17]{BH07}),
and independently Sakai (cf.\ \cite[Theorem 3.7]{Sakai07}).

\blem\label{scof}
The following are equivalent:
\be
\itm $C_p(X,\{0,1\})$ is an $\alpha_1$ space;
\itm For each family $\sseq{\cU_n}$ of pairwise disjoint clopen point-cofinite covers of $X$,
there are cofinite $\cV_n\sbst\cU_n$, $n\in\N$, such that $\bigcup_n\cV_n$ is a point-cofinite cover
of $X$.
\qed
\ee
\elem

A function $f$ with domain $X$ is a \emph{discrete limit} of functions $f_n$, $n\in\N$, if for each
$x\in X$, $f_n(x)=f(x)$ for all but finitely many $n$.

Each bijectively enumerated family $\cU=\{U_n:n\in\N\}$ of subsets of a set $X$
induces a \emph{Marczewski map} $\cU:X\to P(\N)$ defined by
$$\cU(x)=\{n\in\N : x\in U_n\}$$
for each $x\in X$.
The main step in our proof is the following.

\blem\label{alpha1.5marc}
Assume that $C_p(X,\{0,1\})$ is an $\alpha_1$ space, and
$\cU=\sseq{U_n}$ is a bijectively enumerated family of open subsets of $X$.
Then the Marczewski map $\cU:X\to P(\N)$ is a discrete limit of continuous functions.
\elem
\bpf
First, consider the case where for each $n$, $U_n$ is not clopen.

For each $n$, write $U_n$ as a union $\Union_m C^n_m$ of nonempty disjoint clopen sets.
We may assume that the partitions are disjoint:
Inductively, for each $n=2,3,\dots$, consider the elements $C^n_m$, $m\in\N$, of the $n$th
partition. For each $m$, if $C^n_m$ appears in the partition of $U_k$ for some $k<n$,
merge (in the $n$th partition) $C^n_m$ with some other element of the $n$th partition.
Continue in this manner until the $n$th partition is disjoint from all previous partitions.

Thus, the families $\cU_n = \{(C^n_m)\comp : m\in\N\}$ are disjoint clopen point-cofinite covers of $X$.
By Lemma \ref{scof},
there are $k_n$, $n\in\N$, and
subsets $\cV_n=\{(C^n_m)\comp : m\ge k_n\}\sbst\cU_n$, $n\in\N$,
such that $\Union_n\cV_n$ is a point-cofinite cover of $X$. In other words,
$$\cV=\left\{\bigcap^\oo_{m=k_n} (C^n_m)\comp : n\in\N\right\}$$
is a point-cofinite cover of $X$.

For each $n,m$, let
$$U^n_m = \Union_{i=1}^{\max\{m,k_n\}}C^n_i.$$
For each $m$, define $\Psi_m:X\to P(\N)$ by
$$\Psi_m(x)=\{n : x\in U^n_m\}.$$
As each $U^n_m$ is clopen, $\Psi_m$ is continuous.
It remains to prove that, viewed as a Marczewski map, $\cU$ is a discrete limit of the maps $\Psi_m$, $m\in\N$.

Fix $x\in X$.
Let $N$ be such that $x\in \bigcap_{m=k_n}^\oo (C^n_m)\comp$ for all $n\ge N$.
For each $n<N$ with $x\in U_n$, let $m_n$ be such that
$x\in U^n_{m_n}$. Set $M=\max\{m_n : n<N\}$.

Fix $m\ge M$. We show that $x\in U^n_m$ if, and only if, $x\in U_n$.
One direction follows from $U^n_m\sbst U_n$.
To prove the other direction, assume that $x\in U_n$, and consider the two possible
cases:
If $n<N$, then $x\in U^n_m$ because $m\ge M\ge m_n$, and we are done.
Thus, assume that $n\ge N$. Then $x\in \bigcap_{i=k_n}^\oo (C^n_i)\comp$.
As $x\in U_n=\Union_mC^n_m$, it follows that $x\in \bigcup_{i=1}^{k_n-1} C^n_i\sbst U^n_m$.

Thus, for each $x\in X$ there is $M$ such that $\Psi_m(x)=\cU(x)$ for all $m\ge M$.
This completes the proof in the case that no $U_n$ is clopen.

For the remaining case, let $I\sbst\N$ be the set of all $n$ such that $U_n$ is not clopen.
The previous case shows that $\cU_I=\{U_n : n\in I\}$, viewed as a Marczewski function
from $X$ to $P(I)$, is a discrete limit of continuous functions
$\Psi_m:X\to P(I)$.

For each $m$, define $\Phi_m:X\to P(\N)$ by
$$\Phi_m(x) =\{n : (n\in I\mbox{ and }n\in\Psi_m(x))\mbox{ or }(n\nin I\mbox{ and }x\in U_n)\}.$$
Then $\cU$ is a discrete limit of the continuous functions $\Phi_m$, $m\in\N$.
\epf

As $X$ satisfies item (2) of Lemma \ref{scof}, it satisfies $\ufin(\rmO,\Gamma)$:
Refine each given cover to a clopen cover, turn it to a clopen point-cofinite cover by taking
finite unions, and make the point-cofinite covers disjoint.

Assume that $\cU$ is a countable family of open subset of $X$.
By Lemma \ref{alpha1.5marc}, the Marczewski map $\cU:X\to P(\N)$ is a discrete
limit of continuous functions $\Psi_n$.

Clearly, every discrete limit is a quasi-normal limit.
The proof of \cite[Theorem 4.8]{BRR91} actually establishes the following.

\blem\label{P}
Assume that $P$ is a property of topological spaces, which is preserved by taking
closed subsets, continuous images
and countable unions. If $X$ has the property $P$ and
$\Psi:X\to Y$ is a quasi-normal limit of continuous functions
into a metric space $Y$, then $\Psi[X]$ has the property $P$.
\elem
\bpf
Let $\Psi_n$, $n\in\N$, be continuous functions as in the premise of the lemma, and
let $\epsilon_n$, $n\in\N$, be as in the definition of quasi-normal convergence.
For each $k$,
$$X_k = \{x\in X : (\forall n,m\ge k)\ d(\Psi_n(x),\Psi_m(x))\le \epsilon_n+\epsilon_m\}$$
is a closed subset of $X$, and the functions $\Psi_n$ converge to $\Psi$ uniformly
on $X_k$. Thus, $\Psi$ is continuous on $X_k$, and therefore $\Psi[X_k]$ has the
property $P$.

Now, $X=\Union_kX_k$, and therefore $\Psi[X]=\Union_k\Psi[X_k]$ has the property $P$.
\epf

It follows that for each countable family $\cU$ of open subsets of $X$, $\cU[X]$ satisfies $\ufin(\rmO,\Gamma)$.

Let $\rmF,\rmB$ denote the families of all countable closed and all countable Borel covers of $X$,
respectively. Similarly, let $\FG, \BG$ denote the families of all countable closed point-cofinite covers of $X$
and all Borel point-cofinite covers of $X$.
Following is a striking result of Bukovsk\'y, Rec\l{}aw, and Repick\'y  \cite{BRR91}.
In their terminology, it tells that the family of closed subsets of $X$
is weakly distributive if, and only if, the same holds for the family of Borel subsets of $X$.
In the language of selection principles, this result has the following compact form.

\blem[Bukovsk\'y--Rec\l{}aw--Repick\'y \cite{BRR91}]\label{surprise1}
$\ufin(\rmF,\FG)=\ufin(\rmB,\BG)$.
\elem
\bpf
Assume that $X$ satisfies $\ufin(\rmF,\FG)$.
We first show that $X$ is a $\sigma$ space \cite[Theorem 5.2]{BRR91}.

Assume that $G = \bigcap_n U_n$ where for each $n$,
$U_n\supseteq U_{n+1}$ are open subsets of $X$.
Write, for each $n$,
$$U_n = \Union_{m\in\N}C^n_m,$$
where for each $m$, $C^n_m\sbst C^n_{m+1}$ are closed subsets of $X$.
We may assume that the closed cover $\{C^n_m\cup (X\sm U_n):m\in\N\}$ of $X$ has no finite subcover.\footnote{If
there are infinitely many $n$ for which there is some $m_n$ with $C^n_{m_n}=U_n$, then
$G=\bigcap_n C^n_{m_n}$ is closed and we are done. Otherwise, we can ignore
finitely many $n$ and assume that there are no $n,m$ such that $C^n_m=U_n$ contains $G$.}
As $X$ satisfies $\ufin(\rmF,\FG)$ and each given cover is monotone,
there are $m_n$, $n\in\N$, such that $\sseq{C^n_{m_n}\cup (X\sm U_n)}$ is a closed point-cofinite cover of $X$.
For each $k$ define
$$Z_k = \bigcap_{n = k}^\oo C^n_{m_n}.$$
Then each $Z_k$ is a closed subset of $X$, and $G=\Union_k Z_k$ is $F_\sigma$.
This shows that $X$ is a $\sigma$ space.

Now, assume that $\cU_n\in\rmB$, $n\in\N$. Then for each $n$, each element of $\cU_n$
is $F_\sigma$ and can therefore be replaced by countably many closed sets.
Applying $\ufin(\rmF,\FG)$ to the thus modified covers, we obtain a cover in $\FG$.
For each $n$, extend each of the finitely many chosen elements of the $n$th cover to an $F_\sigma$ set from the
original cover $\cU_n$, to obtain an element of $\BG$ chosen in accordance with the definition of
$\ufin(\rmB,\BG)$.\footnote{This argument, in more general form, appears in \cite[Theorem 2.1]{BRR01}.}
\epf

\blem\label{marcz}
The following are equivalent:
\be
\itm $X$ satisfies $\ufin(\rmB,\BG)$;
\itm For each countable family $\cU$ of open subsets of $X$,
$\cU[X]$ satisfies $\ufin(\rmO,\Gamma)$;
\itm For each countable family $\cC$ of closed subsets of $X$,
$\cC[X]$ satisfies $\ufin(\rmO,\Gamma)$.
\ee
\elem
\bpf
$(2\Iff 3)$ Use the auto-homeomorphism of $P(\N)$ defined by mapping a set to its
complement.

$(1\Impl 2)$ The Marczewski map $\cU:X\to P(\N)$ is Borel.
It is easy to see that $\ufin(\rmB,\BG)$ is preserved by Borel images \cite{CBC}.
Thus, $\cU[X]$ satisfies
$\ufin(\rmB,\BG)$, and in particular $\ufin(\rmO,\Gamma)$.

$(3\Impl 1)$ By Lemma \ref{surprise1}, it suffices to show that $X$ satisfies $\ufin(\rmF,\FG)$.
For each $\cC=\sseq{C_n}\in \rmF$ which does not contain a finite subcover, $\sseq{\bigcup_{m\le n}C_m}\in\FG$.
Thus, $\ufin(\rmF,\FG)=\ufin(\FG,\FG)$,\footnote{This statement holds in a more general form \cite{coc2}.}
and we prove the latter property.

Let $\cC_n=\{C^n_m : m\in\N\}$, $n\in\N$, be bijectively enumerated closed point-cofinite covers of $X$
which do not contain finite subcovers. We may assume that these covers are pairwise disjoint \cite{coc1}.

Let $\cC=\Union_n\cC_n$, and consider the Marczewski map $\cC:X\to P(\N\x\N)$ defined by
$$\cC(x)=\{(n,m) : x\in C^n_m\}$$
for all $x\in X$.
For each $(n,m)$, $O_{(n,m)}=\{A\sbst\N\x\N : (n,m)\in A\}$ is an open subset of $P(\N\x\N)$,
and for each $n$, $\cU_n=\{O_{(n,m)}: m\in\N\}$ is an open cover of $\cC[X]$ that
does not contain a finite subcover.
As $\cC[X]$ satisfies $\ufin(\rmO,\Gamma)$, there are $k_n$, $n\in\N$, such that
$\{\Union_{m<k_n} O_{(n,m)} : n\in\N\}$ is a point-cofinite cover of $\cC[X]$.
Then $\{\Union_{m<k_n} C^n_m : n\in\N\}$ is a point-cofinite cover of $X$ (it is infinite because
$X$ does not appear there as an element).\footnote{The
argument is standard: For each finite family of proper subsets of $X$, there is a finite subset
of $X$ not contained in any member of this family.}
\epf

By Lemma \ref{marcz}, $X$ satisfies $\ufin(\rmB,\BG)$.
It remains to observe the following. For the reader's convenience, we
reproduce the proof of the implication needed in the present proof.

\blem[Bartoszy\'nski--Scheepers \cite{BarSch93}]\label{ScTs}
$X$ satisfies $\ufin(\rmB,\BG)$ if, and only if,
each Borel image of $X$ in $\NN$ is bounded.
\elem
\bpf
$(\Impl)$ Assume that $Y\sbst\NN$ is a Borel image of $X$.
Then $Y$ satisfies $\ufin(\rmB,\BG)$. By taking the image of $Y$
under the continuous mapping $f(n)\mapsto f(1)+\dots+f(n)$
defined on $\NN$ , we may assume that all elements in $Y$ are nondecreasing.

We first consider the trivial case:
There is an infinite $I\sbst\N$ such that for each $n\in I$, $F_n=\{f(n) : f\in Y\}$
is finite. For each $n$, let $m\in I$ be minimal such that $n\le m$, and define $g(n)=\max F_m$.
Then $Y$ is bounded by $g$.

Thus, assume that there is $N$ such that for each $n\ge N$, $\{f(n) : f\in Y\}$ is infinite.
For all $n,m$, consider the open set $U^n_m=\{f\in Y : f(n)\le m\}$.
Then for each $n\ge N$, $\cU_n=\{U^n_m : m\in\N\}$ is an open point-cofinite cover of $Y$.
Apply $\ufin(\rmB,\BG)$ to obtain for each $n\ge N$ a finite set $F_n\sbst\N$,
such that $\{\Union_{m\in F_n} U^n_m  : n\in\N\}$ is a point-cofinite cover of $Y$.
Define $g\in\NN$ by $g(n)=\max F_n$ for each $n\ge N$ (and arbitrary for $n<N$).
Then $Y$ is bounded by $g$.
\epf

This completes the proof of Theorem \ref{Borel}.
\epf

\brem
Let $A\as B$ mean that $A\sm B$ is finite. A \emph{semifilter} is a family $\cF$ of infinite
subsets of $\N$ such that for each $A\in\cF$ and each $B\sbst\N$ such that $A\as B$, we have
that $B\in\cF$.
In \cite{SF1} it is proved that if in item (2) of Theorem \ref{marcz} we replace $\cU[X]$ with
the semifilter it generates, then we obtain a characterization of $\ufin(\rmO,\Gamma)$.
Theorem \ref{marcz} shows that moving to the generated semifilter is essential to obtain this
result, since $\ufin(\rmB,\BG)$ is strictly stronger than $\ufin(\rmO,\Gamma)$.
\erem

\section{Applications}

\subsection{QN spaces}
We begin with a straightforward proof of one implication in the following
theorem (which answers in the affirmative Problem 2 of Scheepers \cite{alpha_i}).
Because of the importance of this result, we also supply a proof for
the other implication.

\bthm[Sakai \cite{Sakai07}, Bukovsk\'y--Hale\v{s} \cite{BH07}]\label{SBH}
$X$ is a QN space if, and only if, $C_p(X)$ is an $\alpha_1$ space.
\ethm
\bpf
$(\Leftarrow)$ This is Theorem 4 of \cite{alpha_i}.
Using Theorem \ref{Borel} this becomes straightforward:
Assume that $C_p(X)$ is an $\alpha_1$ space.
Given $f_n$, $n\in\N$, converging pointwise to $0$,
define a Borel function $\Psi:X\to\NN$ by
$$\Psi(x)(n) = \min\{k : (\forall m\ge k)\ |f_m(x)|<1/n\}.$$
By Theorem \ref{Borel}, $\Psi[X]$ is bounded by some $g\in\NN$.
For each $x\in X$ and all but finitely many $n$, $|f_m(x)|<1/n$ for each $m\ge g(n)$.
For each $n$ and each $m$ with $g(n)\le m<g(n+1)$, take $\epsilon_m=1/n$.

$(\Impl)$ Assume that $X$ is a QN space, and $\sseq{A_n}$ is a countable sheaf at $f\in C_p(X)$.
We may assume that $f$ is the zero function, and that the image of each member of each $A_n$
is contained in the unit interval $[0,1]$.

For each $n$, enumerate $A_n=\{f^n_m : m\in\N\}$ bijectively.
For each $m$, define $g_m\in C_p(X)$ by
$$g_m(x) = \sup\{f^n_m(x)/n : n\in\N\}$$
for all $x\in X$.
Then $\{g_m : m\in\N\}$ converges pointwise to the constant zero function.
As $X$ is a QN space, there are positive $\epsilon_m$, $m\in\N$, converging to $0$,
such that $X$ it is the increasing union of the sets
$$X_n=\{x\in X : (\forall m\ge n)\ g_m(x)\le \epsilon_m\}.$$
For each $n$, choose $m_n$ such that $n\epsilon_m\le 1/n$ for all $m\ge m_n$.
We claim that the amalgamation $B=\Union_n\{f^n_m : m\ge m_n\}$ converges
pointwise to the zero function.
Indeed, fix $x\in X$ and a positive $\epsilon$. Take $N$ such that $x$
belong to $X_N$ (and thus to all $X_k$ with $k\ge N$) and such that $1/N\le \epsilon$.
For each $n\ge N$ and each $m\ge m_n$,
$$f^n_m(x)\le n\cdot g_m(x)\le n\epsilon_m\le 1/n\le \epsilon.$$
And for each $n<N$, there are only finitely many $m$ such that $f^n_m(x)>\epsilon$.
Thus, for all but finitely many $f\in B$, $f(x)\le\epsilon$.
\epf

A beautiful direct (but tricky) proof for $(\Leftarrow)$ of Theorem \ref{SBH}
was recently discovered by Bukovsk\'y \cite{BukLower}.

Theorems \ref{Borel} and \ref{SBH} solve in the affirmative Problem 22 from \cite{BH07}.

\bcor\label{QNB}
$X$ is a QN space if, and only if, each Borel image of $X$ in $\NN$ is bounded.\hfill\qed
\ecor

\bthm[Rec\l{}aw \cite{Rec97}]\label{hQN}
The QN property is hereditary.
\ethm
\bpf
The property of having bounded Borel images in $\NN$ is hereditary.
\epf

Answering Question 5.8 of Shakhmatov \cite{Shakhmatov02} (attributed to Scheepers),
Sakai \cite{Sakai07} and independently Bukovsk\'y--Hale\v{s} \cite{BH07}
gave a characterization of the QN property in terms of covering properties of $X$.
Their characterization uses the new \emph{Ko\v{c}inac $\alpha_1$} selection principle \cite{KocAlpha}.
Theorem \ref{ScTs} and Corollary \ref{QNB} give a new
characterization in terms of the classical Hurewicz selection principle: $\ufin(\rmB,\BG)$.
This selection hypothesis can be stated in a more elegant manner. For families of covers
$\scrA,\scrB$ of $X$, define
\begin{description}
\itm[$\sone(\scrA,\scrB)$]
Whenever $\cU_1,\cU_2,\dots\in\scrA$,
there exist elements $U_n\in\cU_n$, $n\in\N$,
such that $\sseq{U_n}\in\scrB$.
\end{description}
Then $\ufin(\rmB,\BG)=\sone(\BG,\BG)$ \cite{CBC}. By Lemma \ref{surprise1}, also
$\ufin(\rmF,\FG)=\sone(\FG,\FG)$. (This can also be proved directly.)
We obtain the following new characterizations.

\bcor\label{sonechar}
The following are equivalent:
\be
\itm $C_p(X)$ is an $\alpha_1$ space;
\itm $X$ is a $QN$ space;
\itm $X$ satisfies $\sone(\FG,\FG)$;
\itm $X$ satisfies $\sone(\BG,\BG)$.\hfill\qed
\ee
\ecor

\subsection{Convergent sequences of Borel functions}
Let $B_p(X)$ be the space of all \emph{Borel} real-valued
functions on $X$, with the topology of pointwise convergence.
We obtain the surprising result, that if $C_p(X)$ is an $\alpha_1$ space,
then so is $B_p(X)$. This is \emph{not} provably the case for \Arh{}'s properties
$\alpha_2$, $\alpha_3$, and $\alpha_4$.
A topological space $Y$ is an \emph{$\alpha_2$ space}
if it satisfies $\sone(\Gamma_y,\Gamma_y)$ for each $y\in Y$, where
$\Gamma_y$ is the family of all sequences converging to $y$.
For the definitions of $\alpha_3$ and $\alpha_4$, see e.g.\ \cite{alpha_i}.

\bcor\label{Bp1}
The following are equivalent:
\be
\itm Each Borel image of $X$ in $\NN$ is bounded;
\itm $C_p(X)$ is an $\alpha_1$ space;
\itm $B_p(X)$ is an $\alpha_1$ space;
\itm $B_p(X)$ is an $\alpha_2$ space;
\itm $B_p(X)$ is an $\alpha_3$ space;
\itm $B_p(X)$ is an $\alpha_4$ space.
\ee
\ecor
\bpf
$(1\Impl 3)$ This is proved verbatim as the proof of $(2\Impl 1)$ in Theorem \ref{Borel}.

$(3\Impl 2)$ is evident.

$(2\Impl 1)$ is due to the mentioned result of Scheepers,
and the equivalence of being a $QN$ space and (1).

$(4\Iff 5\Iff 6)$ is proved as in Gerlits--Nagy \cite{GNFr} or Scheepers' \cite{alpha_i}
(in fact, the Borel case is easier).

$(3\Impl 4)$ is evident.

$(4\Impl 1)$ It suffices to show that $X$ satisfies $\sone(\BG,\BG)$.
Given $\cU_n\in\BG$, $n\in\N$, we have that for each $n$, $A_n = \{\chi_U : U\in\cU_n\}\sbst B_p(X)$
converges pointwise to $0$. Applying $\alpha_2$, let $U_n\in\cU_n$, $n\in\N$, be such that
$\chi_{U_n}$ converges pointwise to $0$. Then $\sseq{U_n}$ is a point-cofinite cover of $X$.
\epf

\subsection{Almost continuous functions}
A function $f:X\to Y$ is \emph{almost continuous} \cite{ArhBok93} if
for each nonempty $A\sbst X$, the restriction of $f$ to $A$
has a point of continuity. $AC_p(X)$ is the space of all almost continuous real valued functions on $X$,
with the topology of pointwise convergence \cite{BokMal00}.

If $X$ and $Y$ are Tychonoff and $f:X\to Y$ is almost continuous,
then for each $A\sbst X$ the set of
points of continuity of the restriction of $f$ to $A$ is open dense in $A$
\cite{BokMal00}.
Each function with the latter property is Borel \cite{Vinokurov85}.
Thus, $C_p(X)\sbst AC_p(X)\sbst B_p(X)$.

\bcor
$AC_p(X)$ is an $\alpha_1$ space if, and only if, $C_p(X)$ is an $\alpha_1$ space.
\ecor

\subsection{wQN spaces and the Scheepers Conjecture}
$X$ is a \emph{wQN space} \cite{BRR91} if each sequence of continuous real-valued functions on $X$
converging pointwise to zero has a subsequence converging to zero quasi-normally.

Two fundamental problems concerning wQN spaces appear in the literature:
In \cite[page 269]{alpha_i}, \cite[Problem 23]{BH07}, and \cite[Problems 10.3--10.4]{Buk07},
we are asked whether, consistently, every wQN space is a QN space.
The \emph{Scheepers Conjecture} \cite{wqn} asserts that $X$ is a wQN space if, and only if,
$X$ satisfies $\sone(\Gamma,\Gamma)$. It is still open whether the Scheepers Conjecture is provable.
A striking result of Dow gives a positive answer to the first problem, and a consistently
positive answer to the second.

\bthm[Dow \cite{Dow90}]\label{Dow}
In the Laver model, each $\alpha_2$ space is an $\alpha_1$ space.
\ethm

Let $\CG$ denote the family of all clopen point-cofinite covers of $X$.
Clearly, $\sone(\Gamma,\Gamma)$ implies $\sone(\CG,\CG)$.

\bcor\label{laver}
In the Laver model:
\be
\itm $\sone(\BG,\BG)=\sone(\CG,\CG)$.
\itm $X$ is a wQN space if, and only if, $X$ is a QN space.
\itm The Scheepers Conjecture holds.
\ee
In particular, these assertions are (simultaneously) consistent.
\ecor
\bpf
(1) Using the correspondence described just before Lemma \ref{scof},
we have that $C_p(X)$ is an $\alpha_2$ space if, and only if, $X$ satisfies $\sone(\CG,\CG)$.
Thus, if $X$ satisfies $\sone(\CG,\CG)$, then by Dow's Theorem \ref{Dow}, $C_p(X)$ is an $\alpha_1$ space.
By Theorem \ref{Borel} $X$ satisfies $\sone(\BG,\BG)$.

(2) Assume that $X$ is a wQN space. Then $C_p(X)$ is an $\alpha_2$ space \cite{CpHure}.
By Dow's Theorem \ref{Dow}, $C_p(X)$ is an $\alpha_1$ space. By Theorem \ref{SBH},
$X$ is a QN space.

(3) $\sone(\Gamma,\Gamma)$ implies (in ZFC) being a wQN space \cite{wqn}.
Now, back in the Laver model, assume that $X$ is a wQN space.
By (2), $X$ is a QN space. By Corollary \ref{sonechar}, $X$ satisfies $\sone(\BG,\BG)$, and
in particular $\sone(\Gamma,\Gamma)$.
\epf

\brem
In \cite{Sakai07, BH07} it is shown that $X$ is a wQN space if, and only if, $X$ satisfies
$\sone(\CG,\CG)$. Using this, (2) and (3) follow immediately from Corollary \ref{laver}(1).
\erem

\subsection{$\overline{\mbox{QN}}$ spaces and $M$ spaces}
$X$ is a \emph{$\overline{\mbox{QN}}$ space} \cite{BRR01} if each real-valued function (not necessarily continuous)
on $X$ which is a pointwise limit of a sequence of continuous functions, is in fact a quasi-normal
limit of those functions.

The following result is immediate from Theorem \ref{Borel} and
\cite[Theorem 5.10(9)]{BRR01}. For completeness, we give a simple, direct proof.

\bthm\label{BorelQN}
The following are equivalent:
\be
\itm $X$ is a $\overline{\mbox{QN}}$ space;
\itm $X$ is a QN space;
\itm Each sequence of Borel functions converging pointwise to $0$, converges to $0$ quasi-normally;
\itm Each sequence of Borel functions converging pointwise to any function, converges quasi-normally to this
function.
\ee
\ethm
\bpf
$(1\Impl 2)$ is immediate.

$(2\Impl 3)$ Assume (2). By Theorem \ref{Borel}, each Borel image of $X$ in $\NN$ is bounded.
Thus, an argument verbatim as in our proof of $(\Leftarrow)$ of Theorem \ref{Borel}
gives (3).

$(3\Impl 4)$ the limit function $f$ is also Borel, and $f_n-f$ converges to $0$.

$(4\Impl 1)$ is immediate.
\epf

This shows that the first assumption in \cite[Theorem 5.10(9)]{BRR01} is not needed.
It also answers \cite[Problem~6.11]{BRR01} in the positive.
Based on \cite[Theorem 6.9]{BRR01} and improving it, we also obtain the following
solution of \cite[Problem 6.10]{BRR01}.

\bcor\label{M}
Every QN space is an $M$ space.\hfill\qed
\ecor

The definition of $M$ space is available at \cite{BRR01}.


\subsection{wQN$_*$ spaces}
A space $X$ is \emph{wQN$_*$} if each sequence of lower semi-continuous real-valued functions on $X$
converging pointwise to zero has a subsequence converging to zero quasi-normally.
In his talk at the \emph{Third Workshop on Coverings, Selections, and Games in Topology} (Serbia, April 2007),
Bukovsk\'y defined wQN$_*$ spaces and described his recent investigations of this property and its
upper semi-continuous variant. The main problem he posed was: Is every QN space a wQN$_*$ space?

\bthm\label{lower}
Every QN space is a wQN$_*$ space.
\ethm
\bpf
Every lower semi-continuous function is Borel. Use Theorem \ref{BorelQN}.
\epf

Bukovsk\'y has later proved that the converse implication in Theorem \ref{lower} also holds \cite{BukLower},
and therefore the notions coincide (with one another and with having bounded Borel images).

\subsection{Bounded-ideal convergence spaces}
The notion of ideal convergence originates in works of Steinhaus and Fast on statistical convergence,
and was generalized by Bernstein, Kat\v{e}tov, and others (see \cite{FMRS07} for an introduction).
The following definitions are as in Jasinski--Rec\l{}aw \cite{JaRec06}.
Let $D$ be a countable set, and $\cI\sbst P(D)$ be an ideal (i.e., $\cI$ contains all singletons and
is closed under taking subsets and finite unions). $\cI^*$ denotes the filter $\{D\sm A : A\in\cI\}$
dual to $\cI$. A sequence $\{r_d\}_{d\in D}$ of real numbers \emph{$\cI$-converges} to $0$ if
for each positive $\epsilon$, $\{d\in D : |r_d|<\epsilon\}\in \cI^*$.
A sequence $\{f_d\}_{d\in D}$ of continuous real-valued functions on $X$ $\cI$-converges to
$0$ if for each $x\in X$, the sequence of real numbers $\{f_d(x)\}_{d\in D}$ $\cI$-converges to $0$.
A space $X$ has the \emph{$\cI$-convergence property} if for each sequence $\{f_d\}_{d\in D}$ of
continuous real-valued functions on $X$ which $\cI$-converges to $0$, there is $A\in\cI^*$
such that $\{f_d\}_{d\in A}$ converges pointwise to $0$.

We will use the following.

\blem\label{distinct}
In the definition of the $\cI$-convergence property, it suffices to consider only
sequences of \emph{distinct} elements.
\elem
\bpf
Let $\{f_d\}_{d\in D}$ be given.
Enumerate $D=\sseq{d_n}$ bijectively. For each $n$,
as the functions $f_{d_n}+1/m$, $m\in \N$, are all distinct, there is
$m(d_n)\in \N$ such that $m(d_n)\ge n$ and $f_{d_n}+1/m(d_n)\nin\{f_{d_1}+1/m(d_1),\dots,f_{d_{n-1}}+1/m(d_{n-1})\}$.

It is easy to see that $\{f_d\}_{d\in D}$ $\cI$-converges to $0$ if, and only if,
$\{f_d+1/m(d)\}_{d\in D}$ $\cI$-converges to $0$.
\epf

We use these definitions for $D=\N\x\N$.
For $h\in\NN$, define $A_h=\{(n,m) : m\le h(n)\}$. The family $\{A_h :h\in\NN\}$ is closed
under finite intersections, and generates the \emph{bounded-ideal}
$$\cI_b=\{B\sbst\N\x\N : (\exists h\in\NN)\ B\sbst A_h\}.$$
$X$ has the \emph{bounded-ideal convergence property} if it has the $\cI_b$-convergence
property.

The bounded-ideal, which is also the Fubini product $\emptyset\x\mathrm{Fin}$ of the
trivial ideal and the ideal of finite sets, plays a central role in studies
of ideal convergence. For each analytic $P$-ideal $\cI$, if any $X\sbst\R$
not having Lebesgue measure zero has the $\cI$ ideal convergence, then
$\cI$ is isomorphic to $\cI_b$ \cite{JaRecCM}. For additional uses of this
ideal and its associated convergence, see \cite{FMRS07}.

Jasinski and Rec\l{}aw \cite{JaRec06} proved that every Sirepi\'nski set has the bounded-ideal convergence
property, and that if $X$ has the bounded-ideal convergence property, then $X$ is a $\sigma$ space.
Both of these assertions follow at once from the following.

\bthm\label{Ib}
The following are equivalent:
\be
\itm $X$ has the bounded-ideal convergence property;
\itm $C_p(X)$ is an $\alpha_1$ space;
\itm Each Borel image of $X$ in $\NN$ is bounded.
\ee
\ethm
\bpf
By Theorem \ref{Borel}, it suffices to show that $(1\Iff 2)$.

$(2\Impl 1)$ Assume that $\{f_{(n,m)}\}_{(n,m)\in\N\x\N}$ $\cI_b$-converges to $0$.
By Lemma \ref{distinct}, we may assume that the elements $f_{(n,m)}$, $(n,m)\in\N\x\N$,
are distinct.

For each $x\in X$ and each positive $\epsilon$,
$\{(n,m) : |f_{(n,m)}(x)|<\epsilon\}\in\cI_b^*$, that is,
there is $h\in\NN$ such that $\{(n,m) : |f_{(n,m)}(x)|<\epsilon\}\spst (\N\x\N)\sm A_h$.
Thus, $|f_{(n,m)}(x)|<\epsilon$ for all $n,m\in\N$ such that $h(n)<m$.
It follows that for each $n$, $\{f_{(n,m)}\}_{m\in\N}$ converges pointwise to $0$.

As $C_p(X)$ is an $\alpha_1$ space, there is for each $n$ a number $h(n)\in\N$ such
that $\{f_{(n,m)} : n,m\in\N, m>h(n)\}$ converges pointwise to $0$, and
since its enumeration is bijective, the sequence $\{f_{(n,m)}\}_{(n,m)\in (\N\x\N)\sm A_h}$
also converges pointwise to $0$.
As $(\N\x\N)\sm A_h\in\cI_b^*$, this shows that $X$ has the bounded-ideal convergence property.

$(1\Impl 2)$ Assume that for each $n$ the sequence $\{f_{(n,m)}\}_{m\in\N}$ converges pointwise to $0$.
For each $x\in X$, each positive $\epsilon$, and each $n$, there is $h(n)\in\N$
such that $|f_{(n,m)}(x)|<\epsilon$ for all $m>h(n)$.
Thus, $\{(n,m) : |f_{(n,m)}(x)|<\epsilon\}\spst (\N\x\N)\sm A_h$,
that is, $\{f_{(n,m)}\}_{(n,m)\in\N\x\N}$ $\cI_b$-converges to $0$.

By the bounded-ideal convergence property, there is $h\in\NN$ such that
$\{f_{(n,m)}\}_{(n,m)\in (\N\x\N)\sm A_h}$ converges pointwise to $0$,
and therefore so does the sheaf amalgamation $\{f_{(n,m)} : n,m\in\N, m>h(n)\}$ (which can
be enumerated as a subsequence of $\{f_{(n,m)}\}_{(n,m)\in (\N\x\N)\sm A_h}$.
\epf

\begin{cor}
The following are equivalent:
\be
\itm $X$ has the bounded-ideal convergence property;
\itm For each sequence $\{f_d\}_{d\in \N\x\N}$ of
Borel real-valued functions on $X$ which $\cI_b$-converges to $0$, there is $A\in\cI_b^*$
such that $\{f_d\}_{d\in A}$ converges pointwise to $0$;
\itm For each sequence $\{f_d\}_{d\in \N\x\N}$ of
Borel real-valued functions on $X$ which $\cI_b$-converges to a Borel function $f$, there is $A\in\cI_b^*$
such that $\{f_d\}_{d\in A}$ converges pointwise to $f$;
\itm For each sequence $\{f_d\}_{d\in \N\x\N}$ of
Borel real-valued functions on $X$ which $\cI_b$-converges to a function $f$, there is $A\in\cI_b^*$
such that $\{f_d\}_{d\in A}$ converges pointwise to $f$.
\ee
\end{cor}
\bpf
$(1\Impl 2)$ Replace ``continuous'' by ``Borel'' in the proof of Theorem \ref{Ib} and use Theorem \ref{Bp1}.

$(2\Impl 3)$ $B_p(X)$ is a topological group, and in particular homogeneous.

$(3\Impl 4)$ The assumption in (4) implies, in particular, that $f$ is a pointwise limit of
$\{f_{(1,m)}\}_{m\in\N}$. Thus, $f$ is Borel.
\epf

\subsection{Bounded Baire-class $\alpha$ images}
Continuous functions and Borel functions are the extremal notions in the \emph{Baire hierarchy} of functions:
A real-valued function $f$ is of \emph{Baire-class $0$} if it is continuous. For $0<\alpha\le \aleph_1$,
$f$ is of \emph{Baire-class $\alpha$} if $f$ is the pointwise limit of a sequence of functions, each
of Baire-class smaller than $\alpha$.
$f$ is Borel if, and only if, $f$ is of Baire class $\aleph_1$ (see \cite{Kechris}).
A natural question in light of our study is: Which spaces $X$ have the property that
each Baire-class $\alpha$ image of $X$ in $\NN$ is bounded?

\bthm
For each $\alpha>0$, the following are equivalent:
\be
\itm Each Baire-class $\alpha$ image of $X$ in $\NN$ is bounded;
\itm Each Borel image of $X$ in $\NN$ is bounded.
\ee
\ethm
\bpf
Assume that each Baire-class $1$ image of $X$ in $\NN$ is bounded.
Baire-class $1$ functions are exactly the $F_\sigma$-measurable functions.

One way to proceed is using Lemma \ref{marcz}, since
for each bijectively enumerated family of closed sets $\cC=\sseq{C_n}$,
the corresponding Marczewski function is $F_\sigma$-measurable (and by the
proof of Lemma \ref{marcz}, we may assume that for each $x\in X$, $\cC(x)$ is infinite).

However, there is a more direct proof. By Lemma \ref{surprise1}, it suffices to
prove that $X$ satisfies $\ufin(\rmF,\FG)$. Assume that $\cU_n=\{C^n_m : m\in\N\}$,
$n\in\N$, are closed covers of $X$ not containing a finite subcover.
Define $\Psi:X\to\NN$ by
$$\Psi(x)(n)=\min\{m : x\in C^n_m\}$$
for all $n\in\N$. Each basic open subset of $\NN$ is an intersection of finitely many
sets of the form $O^n_m = \{f\in\NN : f(n)=m\}$. As
$\Psi\inv[O^n_m]=C^n_m\sm \Union_{k<m} C^n_k$ is an $F_\sigma$ set for all
$n$ and $m$, $\Psi$ is $F_\sigma$-measurable. Thus, $\Psi[X]$ is bounded by some $g\in\NN$.
Then $\sseq{\Union_{m\le g(n)}C^n_m}$ is a point-cofinite cover of $X$.
\epf

\section{Closing the circle: Continuous bounded images again}

The proof of Theorem \ref{Borel} gives us the following analogue of Theorem \ref{main}.
Say that a set $Y\sbst\NcompactN$ is \emph{bounded} if there is $g\in\NN$
such that for each $f\in Y$ and all but finitely many $n$, $f(n)<\oo$ implies $f(n)\le g(n)$.
This generalizes the standard notions of boundedness in $\NN$ or $\EF$.

\bthm\label{contBorel}
The following are equivalent:
\be
\itm Each Borel image of $X$ in $\NN$ is bounded;
\itm Each continuous image of $X$ in $\NcompactN$ is bounded.
\ee
\ethm
\bpf
$(1\Impl 2)$ Assume that $\Psi:X\to\NcompactN$ is continuous.
Define $d:\NcompactN\to\NN$ by $d(x)(n)=x(n)$ if $x(n)<\oo$, and $d(x)(n)=1$
if $x(n)=\oo$. Then $d\circ\Psi:X\to\NN$ is Borel, and therefore
$d[\Psi[X]]$ is a bounded subset of $\NN$.
Thus, $\Psi[X]$ is a bounded subset of $\NcompactN$.

$(2\Impl 1)$ Assume that each continuous image of $X$ in $\NcompactN$ is bounded.
We first prove that for each bijectively enumerated family $\cU=\sseq{U_n}$ of open sets,
$\cU$ is a discrete limit of continuous functions. The proof is similar to the
proof of Lemma \ref{alpha1.5marc}.
As shown at the end of the proof of Lemma \ref{alpha1.5marc},
we may assume that no $U_n$ is clopen.

For each $n$, write $U_n=\Union_m C^n_m$ as a union of disjoint clopen sets.
Define $\Psi:X\to\NcompactN$ by
$$\Psi(x)(n)=
\begin{cases}
m & x\in U^n_m\\
\oo & x\nin U_n
\end{cases}
$$
Then $\Psi$ is continuous. Let $g\in\NN$ bound $\Psi[X]$.
For each $n,m$, let
$$U^n_m = \Union_{i=1}^{\max\{m,g(n)\}}C^n_i.$$
For each $m$, define a continuous function $\Psi_m:X\to P(\N)$ by
$$\Psi_m(x)=\{n : x\in U^n_m\}.$$
We claim that $\cU$ is a discrete limit of the maps $\Psi_m$, $m\in\N$.

Fix $x\in X$.
Let $N$ be such that for all $n\ge N$, $\Psi(x)(n)<\oo$ implies $\Psi(x)(n)\le g(n)$.
For each $n<N$ with $x\in U_n$, let $m_n$ be such that
$x\in U^n_{m_n}$. Set $M=\max\{m_n : n<N\}$.

Fix $m\ge M$. We show that $n\in\Psi_m(x)$ if, and only if, $x\in U_n$.
One direction follows from $U^n_m\sbst U_n$.
To prove the other direction, assume that $x\in U_n$, and consider the two possible
cases:
If $n<N$, then $x\in U^n_m$ because $m\ge M\ge m_n$, and we are done.
Thus, assume that $n\ge N$. Then $\Psi(x)(n)\le g(n)$, and therefore
$x\in \bigcup_{i=1}^{g(n)} C^n_i\sbst U^n_m$.

Thus, for each $x\in X$ there is $M$ such that for all $m\ge M$, $\Psi_m(x)=\cU(x)$.

Now, each continuous image of $X$ in $\NN$ is bounded because $\NN$ is a subspace of
$\NcompactN$. By the Hurewicz Theorem \ref{hure}, $X$ satisfies $\ufin(\rmO,\Gamma)$,
and by Lemma \ref{P}, so does $\cU[X]$.
By Lemma \ref{marcz}, each Borel image of $X$ in $\NN$ is bounded.
\epf

We therefore obtain the aesthetically pleasing result, that the chain of properties
$$\ufin(\rmB,\BG)\quad\Longrightarrow\quad \mbox{hereditarily-}\ufin(\rmO,\Gamma)\quad\Longrightarrow
\quad \ufin(\rmO,\Gamma)$$
is obtained by requiring bounded continuous images in the chain of subspaces
$$\NcompactN\quad\spst\quad \EF\quad\spst\quad\NN,$$
respectively.

\subsection*{Acknowledgments}
We thank Masami Sakai for his Theorem \ref{sig'}.
We also thank Lev Bukovsk\'y for useful comments and for
making his work \cite{BukLower} available to us prior to its publication,
and the referee for his useful comments.
A part of the present paper was written when the first author was
visiting the Department of Mathematics at the University of Warsaw.
We thank Tomasz Weiss for his kind hospitality during that stimulating period.

\subsubsection*{Note} After the present paper was accepted for publication,
Bukovsk\'y and \v{S}upina devised an alternative, more analytic and less combinatorial,
proof of our main Theorem \ref{Borel} \cite{BukSup12}.

\ed